\documentclass[11pt]{amsart}

\usepackage{amsmath,amssymb}

\makeatletter
\renewcommand{\mod}[1]{\allowbreak \if@display \mkern 12mu \else \mkern 6mu\fi {\operator@font mod}\,\,#1}
\makeatother

\newtheorem{theorem}{Theorem}[section]
\newtheorem{add}[theorem]{Addendum}
\newtheorem{alg}[theorem]{Algorithm}
\newtheorem{claim}[theorem]{Claim}
\newtheorem{conj}[theorem]{Conjecture}
\newtheorem{cor}[theorem]{Corollary}
\newtheorem{dfn}[theorem]{Definition}
\newtheorem{proposition}[theorem]{Proposition}
\newtheorem{lemma}[theorem]{Lemma}
\newtheorem{caution}[theorem]{Caution}
\newtheorem{exa}[theorem]{Example}
\newtheorem{exc}[theorem]{Exercise}
\newtheorem{examples}[theorem]{Examples}
\newtheorem{remark}[theorem]{Remark}
\numberwithin{equation}{section}

\newcommand{\I}{\mathrm{I}}
\newcommand{\II}{\mathrm{II}}
\newcommand{\III}{\mathrm{III}}
\newcommand{\IV}{\mathrm{IV}}

\newcommand{\rd}[1]{\left\lfloor#1\right\rfloor}
\newcommand{\ru}[1]{\left\lceil#1\right\rceil}
\newcommand{\Span}[1]{\left<#1\right>}
\newcommand{\orb}{{}_{\mathrm{orb}}}
\newcommand{\per}{{}_{\mathrm{per}}}
\newcommand{\grow}{{}_{\mathrm{grow}}}
\newcommand{\tensor}{\otimes}
\newcommand{\ot}{\leftarrow}
\newcommand{\bij}{\leftrightarrow}
\newcommand{\onto}{\twoheadrightarrow}
\newcommand{\into}{\hookrightarrow}
\newcommand{\iso}{\cong}

\newcommand{\Aff}{\mathbb A}
\newcommand{\C}{\mathbb C}
\newcommand{\PP}{\mathbb P}
\newcommand{\Gm}{\mathbb G_{\mathrm m}}
\newcommand{\Q}{\mathbb Q}
\newcommand{\Z}{\mathbb Z}
\newcommand{\Ga}{\Gamma}
\newcommand{\De}{\Delta}

\newcommand{\al}{\alpha}
\newcommand{\be}{\beta}
\newcommand{\ga}{\gamma}
\newcommand{\de}{\delta}
\newcommand{\la}{\lambda}
\newcommand{\si}{\sigma}
\newcommand{\Si}{\Sigma}
\newcommand{\ep}{\varepsilon}
\newcommand{\om}{\omega}
\newcommand{\bmu}{\boldsymbol\mu}
 \newcommand{\sia}{\stackrel{a}{\si}}
 \newcommand{\sib}{\stackrel{b}{\si}}

\newcommand{\Oh}{\mathcal O}
\newcommand{\sB}{\mathcal B}
\newcommand{\sC}{\mathcal C}
\newcommand{\sK}{\mathcal K}
\newcommand{\sL}{\mathcal L}
\newcommand{\sN}{\mathcal N}
\newcommand{\sP}{\mathcal P}
\newcommand{\sQ}{\mathcal Q}

\newcommand{\tP}{\widetilde P}

\DeclareMathOperator{\hcf}{hcf}
\DeclareMathOperator{\bhcf}{\mathbf{hcf}}
\DeclareMathOperator{\ch}{ch}
\DeclareMathOperator{\InvMod}{InvMod}
\DeclareMathOperator{\Div}{Div}
\DeclareMathOperator{\Hom}{Hom}
\DeclareMathOperator{\Num}{Num}
\DeclareMathOperator{\Proj}{Proj}
\DeclareMathOperator{\Sing}{Sing}
\DeclareMathOperator{\Spec}{Spec}
\DeclareMathOperator{\RR}{RR}
\DeclareMathOperator{\Td}{Td}
\DeclareMathOperator{\wt}{wt}

\begin{document}

\title{Ice cream and orbifold Riemann--Roch}

\author{Anita Buckley \and Miles Reid \and Shengtian Zhou}

\thanks{Partially funded by Korean Government WCU Grant
R33-2008-000-10101-0. Shengtian Zhou was supported by a University of
Warwick Postgraduate Research Studentship.}
\dedicatory{To Professor Igor Rostislavovich Shafarevich on his 90th birthday}

\begin{abstract}
We give an orbifold Riemann--Roch formula in closed form for the Hilbert
series of a quasismooth polarized $n$-fold $(X,D)$, under the assumption
that $X$ is projectively Gorenstein with only isolated orbi\-fold
points. Our formula is a sum of parts each of which is integral and
Gorenstein symmetric of the same canonical weight; the orbifold parts
are called {\em ice cream functions}. This form of the Hilbert series is
particularly useful for computer algebra, and we illustrate it on
examples of K3 surfaces and Calabi--Yau 3-folds.

These results apply also with higher dimensional orbifold strata (see
\cite{BuSz} and \cite{zh}), although the correct statements are
considerably trickier. We expect to return to this in future
publications.
\end{abstract}



\maketitle

\tableofcontents

\section{Introduction} \label{s!intro}

Reid \cite{YPG} introduced Riemann--Roch (RR) formulas for polarized
orbi\-folds $(X,D)$ with isolated orbifold locus, of the form
\begin{equation}
\label{eq!YPG}
\chi(X,\Oh_X(D))=\RR(X,D) +\sum_{P\in\sB} c_P(D),
\end{equation}
where $\RR(X,D)$ is a Riemann--Roch like expression and the $c_P(D)$
are certain fractional contributions from the orbifold points $\sB$,
depending only on the local type of $(X,D)$. The orbifold RR formula
of \cite{YPG} has found numerous subsequent extensions and
applications; see for example Iano-Fletcher \cite{fl2},
Brown, Alt{\i}nok and Reid \cite{Si}, Buckley and Szendr{\H o}i
\cite{BuSz}, Chen, Chen and Chen \cite{CCC} and Kawakita
\cite{Kawakita}, and we expect these ideas to be equally applicable in
the study of higher dimensional varieties.

A general RR formula for abstract orbifolds was first proved by Kawasaki
\cite{kawasaki} by analytic tools. To{\"e}n \cite{toen} gave another
proof using the algebraic methods of Deligne--Mumford stacks. However,
at present, how to use these abstract results in practice to compute the
dimension of RR spaces is not well understood. To{\"e}n's result was
applied to weighted projective spaces by Nironi \cite{nironi}, to
quasismooth varieties in weighted projective spaces by Zhou \cite{zh}
and to twisted curves by Abramovich and Vistoli \cite{abrvist}. Edidin's
recent treatment \cite{Edi} clarifies orbifold RR considerably; our
results and Zhou's thesis \cite{zh} provide many practical exercises.
Our proof, like that of \cite{YPG}, is based on a reduction to
Atiyah--Singer and Atiyah--Segal equivariant Riemann--Roch
\cite{atiyahseg}, \cite{atiyah}.

Let $D$ be an ample $\Q$-Cartier divisor on a normal projective $n$-fold
$X$ (we usually work over $\C$). The finite dimensional vector spaces
$H^0(X,\Oh_X(mD))$ fit together as a finitely generated graded ring
\begin{equation}\label{eq!RX}
R(X,D)=\bigoplus_{m \ge0} H^0(X,\Oh_X(mD)),
\end{equation}
with $X\iso\Proj R(X,D)$ and the divisorial sheaf $\Oh_X(mD)$ equal to
the character sheaf $\Oh_X(m)$ of the $\Proj$. A surjection from a
graded polynomial ring
\begin{equation}
k[x_0,\dots,x_N] \onto R(X,D)
\quad\hbox{with variables $x_i$ of weight $a_i$}
\end{equation}
corresponds to an embedding
\begin{equation}
i\colon X\iso\Proj R(X,D)\into \PP(a_0,\dots,a_N)
\end{equation}
of $X$ into a weighted projective space as a projectively normal
subscheme.

The Hilbert function $m\mapsto P_m(X,D)=h^0(X,\Oh_X(mD))$ and the {\em
Hilbert series} $P_X(t)=\sum_{m\ge0} P_mt^m$ encode the numerical data
of $R(X,D)$. It is a standard result that $\prod(1-t^{a_i})\cdot P_X(t)$
is a polynomial where, as above, the $a_i$ are the weights of the
generators. The multiplicative group $\Gm$ ($=\C^\times$ if the ground
field is $\C$) has a standard action on the graded ring
$R(X,D)=\bigoplus_{m\ge0} R_m$, with $\la\in\C^\times$ multiplying $R_m$
by $\la^m$. Our aim is a {\em character formula} expressing the Hilbert
series of $R$ in closed form.

\subsection{Plan of the paper} \label{ss!plan}
Chapter~\ref{s!intro} recalls notation and background results from the
literature, and states our Main Theorem~\ref{th!main}.
Chapter~\ref{s!ice} defines the ice cream functions
$P\orb\bigl(\tfrac1r(a_1,\dots,a_n),k_X\bigr)$ as inverse polynomials
modulo $1+t+\cdots+t^{r-1}$ that contain the same information as
Dedekind sums (see especially \ref{ss!Porb}). Chapter~\ref{s!RRD} deals
with the existence of the RR formula for $n$-folds with isolated
orbifold points and the precise nature of the term $\RR(D)$, as a
preliminary to the proof of the main theorem in Section~\ref{s!mpf}.
Section~\ref{K3} relates the new viewpoint of this paper to traditional
formulas for the Hilbert series of K3 surfaces, Fano 3-folds and
canonical 3-folds.

Although this paper mostly deals in isolated orbifold points, our
ultimate aspiration is to find closed expressions for the Hilbert series
of arbitrary orbifolds, having a stratification by orbifold loci of any
dimension. Chapter~\ref{s!hd} discusses briefly what we hope to do in
this direction, and the difficulties associated with positive
dimensional orbifold loci, especially their dissident strata (where the
inertia group jumps); we exemplify this with Buckley's results on
orbifold RR for polarized Calabi--Yau 3-folds \cite{BuSz}. The paper is
backed up by a website
\begin{quote}
 \verb!http://warwick.ac.uk/staff/Miles.Reid/Ice!
\end{quote}
containing additional material that does not fit in the paper, including
implementations of our main algorithms in computer algebra and links to
other papers.

A long term motivating question for us is the ``exact plurigenus
formula'' for Fano 4-folds and canonical 4-folds. We believe that
terminal singularities are intractable in dimension $\ge4$, so it is
unreasonable to hope for a reduction in the style of \cite{YPG} of
general terminal singularities to orbifolds points. Nevertheless, the
orbifold examples provide rich experimental material, and some
modifications of the ideas of this paper should apply more generally.

Our term {\em dissident point} was originally coined in the 1970s as a
reference to Professor Igor Rostislavovich Shafarevich, who created the
Moscow school of algebraic geometry, and who has taught us so much. It
is a pleasure to offer this paper to him as a birthday tribute.

\subsection{Definitions and notation} \label{ss!def}
A {\em Weil divisor} on a normal variety $X$ is a formal linear
combination of prime divisors with integer coefficients. A Weil divisor
$D$ is $\Q$-{\em Cartier} if $mD$ is Cartier for some integer $m>0$.

We write $\bmu_r\subset\Gm$ for the multiplicative group of $r$th roots
of unity, or the cyclic subgroup of $\C^\times$ generated by
$\exp\frac{2\pi i}r$. A {\em cyclic orbifold point} or {\em cyclic
quotient singularity} of type $\frac1r(a_1,\dots,a_n)$ is the quotient
$\pi\colon \Aff^n\to\Aff^n/\bmu_r$, where $\bmu_r$ acts on $\Aff^n$ by
\begin{equation}
\bmu_r\ni\ep \colon (x_1,\dots,x_n)\mapsto (\ep^{a_1}x_1,\dots,\ep^{a_n}x_n).
\end{equation}
We usually assume that no factor of $r$ divides all the $a_i$, which is
equivalent to the $\bmu_r$ action being effective; the orbifold point is
isolated if and only if all the $a_i$ are coprime to $r$. The sheaf
$\pi_*\Oh_{\Aff^n}$ decomposes as a direct sum of divisorial eigensheaves
\begin{equation}
\label{eq!eig}
\sL_i=\bigl\{f \bigm| \ep(f)=\ep^i\cdot f
\hbox{ for all $\ep\in\bmu_r$}\bigr\}
\quad\hbox{for $i\in\Z/r=\Hom(\bmu_r,\Gm)$.}
\end{equation}

The notation $\frac1r(a_1,\dots,a_n)$ refers to {\em polarized} orbifold
points. The orbi\-nates $x_j$ of degree $a_j$ modulo $r$ are local
sections of $\Oh_X(a_j)$, which is locally isomorphic to $\sL_{-a_j}$.
In the terminology of \cite[Definition~8.3]{YPG}, $\Oh_X(1)=\Oh_X(D)$ is
of type ${}_{r-1}\bigl(\frac1r(a_1,\dots,a_n)\bigr)$.

A polarized variety $(X,D)$ is {\em quasismooth} if the corresponding
affine cone $\sC_X=\Spec R(X,D)$ is nonsingular outside the origin. In
this case, the orbi\-fold points of $X$ arise from the orbits of the
group action that are pointwise fixed by a nontrivial isotropy group,
necessarily the cyclic subgroup $\bmu_r\subset\Gm$ for some $r$. In
terms of $(X,D)$, quasismooth holds if and only if $X$ has locally
cyclic quotient singularities $\frac1r(a_1,\dots,a_n)$ and the given
Weil divisor $D=\Oh_X(1)$ generates the local class group
$\Z/r=\Hom(\bmu_r,\Gm)$. Then the local index one cyclic cover defined
by a local identification $\Oh_X(rD)\iso\Oh_X$ is nonsingular.

All our concrete examples are subvarieties in weighted projective
spaces; see Iano-Fletcher \cite{fl2} for definitions and properties. Our
quasismooth assumption implies that $X$ has no orbifold behaviour in
codimension 0 or~1, or is {\em well formed} in the terminology of
\cite{fl2}. This is right here because we work with $n$-folds for
$n\ge2$ with isolated orbifold locus; it means that the orbifold $X$ as
a scheme already knows its orbifold structure, the local universal cover
of $X\setminus\Sing X$. This simplifies the treatment, allowing us to
circumvent the language of stacks and the graded structure sheaf
$\bigoplus_{i\in\Z}\Oh_X(i)$ (cf.\ Canonaco \cite{Ca}). Some of our
examples involve fractional divisors on curves, and we leave the
elementary treatment of the graded structure sheaf
$\bigoplus_{i\in\Z}\Oh_X(i)$ in this case to the conscientious
reader.\footnote{See for example Exercise~\ref{exc!57}. Compare also
Demazure \cite{De} and Watanabe \cite{Wa}; the latter also treats the
graded dualizing sheaf for fractional divisors.}

A polarized variety $(X,D)$ is {\em projectively Gorenstein} if its
affine cone or the corresponding graded ring $R(X,D)$ is Gorenstein. In
this case $\om_X\iso\Oh_X(k_XD)$ for some $k_X\in\Z$, called the {\em
canonical weight} of $(X,D)$, and $H^j(X,\Oh_X(mD))=0$ for all $0<j<\dim
X$ and all $m$. Bruns and Herzog \cite[Corollary 4.3.8]{BrHer} give the
following elementary result:

\begin{lemma}\label{brrherr}
Let $R$ be a graded Gorenstein ring of dimension $\dim R=n+1$ and
canonical weight $k_R$, so that the canonical module of $R$ is
$\om_R=R(k_R)$. Then the Hilbert series $P_R(t)$ of $R$ satisfies
the functional equation
\begin{equation}
t^{k_R}P(\tfrac1t)=(-1)^{n+1}P(t).
\label{eq!Gsym}
\end{equation}
\end{lemma}

We refer to property \eqref{eq!Gsym} of a rational function as {\em
Gorenstein symmetry}. A palindromic polynomial or Laurent polynomial is
Gorenstein symmetric. Examples: $t$ and $t^{-1}+1+t^2+t^3$ are both
palindromic of degree $2$.

\smallskip
\paragraph{\bf Proof} This follows from duality: $R$ is a quotient of a
weighted polynomial ring $A=k[x_0,\dots,x_N]$ with $\wt x_i=a_i$. A
minimal free resolution
\begin{equation}
R \ot F_0 \ot F_1 \ot\cdots\ot F_{\mathrm{cod}}\ot 0,
\end{equation}
has length equal to the codimension $\mathrm{cod}=N-n$, and
$F_{\mathrm{cod}}=A(-\al)$ is the free module of rank one and degree
$-\al$, where $\al=k_R+\sum a_i$ is the {\em adjunction number} for
$X=\Proj R \subset\PP(a_0,\dots,a_N)$. Duality gives
$F_{\mathrm{cod}-i}\iso\Hom_A(F_i,F_{\mathrm{cod}})$ so that, over the
denominator $\prod(1-t^{a_i})$ corresponding to $A=k[x_0,\dots,x_N]$,
the numerator of the Hilbert series is a sum of terms
$t^d+(-1)^{\mathrm{cod}} t^{\al-d}$.
\qed\par\medskip

For quasismooth $X$, the statement corresponds to Serre duality.
However, the proof only uses the definition and basic properties of
Gorenstein graded rings, without further assumptions on the
singularities of $\Spec R$ or $\Proj R$.

Following Mukai \cite{mukai}, we write $c=k_X+n+1$ for the {\em coindex}
of $(X,D)$. By the adjunction formula, the coindex is invariant under
passing to a hyper\-plane section of degree~1. For nonsingular varieties,
we have:
\begin{exa} \label{ex!cix} \rm
\begin{equation}
\begin{tabular}{rc}
projective space $\PP^n$ & has coindex 0; \\[6pt]
a quadric $Q \subset \PP^{n+1}$ & has coindex 1;\\[6pt]
an elliptic curve,
del Pezzo surface \\
or Fano 3-fold of index 2 & has coindex 2; \\[6pt]
a canonical curve,
K3 surface \\
or anticanonical Fano 3-fold & has coindex 3; \\[6pt]
a canonical surface, Calabi--Yau 3-fold \\
or anticanonical Fano 4-fold & has coindex 4.
\end{tabular}
\notag
\end{equation}
\end{exa}

\subsection{The main result} \label{ss!mthm}
For a quasi\-smooth projectively Gorenstein orbifold $(X,D)$ with
isolated orbifold points, Theorem~\ref{th!main} writes the Hilbert
series of $(X,D)$ as a sum of parts, each of which is integral and
Gorenstein symmetric of the same degree $k_X$. We call the orbifold
contributions $P\orb(Q,k_X)$ {\em ice cream functions}; see
\eqref{PorbInteg!eq}. The result expresses $P_X(t)$ in a closed form
that can be calculated readily as a few lines of computer algebra.

\begin{theorem} \label{th!main}
Let $(X,D)$ be a quasismooth orbifold of dimension $n\ge2$. Suppose that
$(X,D)$ is projectively Gorenstein of canonical weight $k_X$, and has
isolated orbifold points
\[
\sB=\left\{Q \hbox{ of type } \tfrac1r(a_1,\dots,a_n)\right\}.
\]

Then the Hilbert series of $X$ is
\begin{equation}
P_X(t)=P_I(t)+\sum_{Q\in\sB} P\orb(Q,k_X)(t),
\end{equation}
where
\begin{enumerate}
\renewcommand{\labelenumi}{(\roman{enumi})}

\item the {\em initial part} has the form
$P_I=\frac{A(t)}{(1-t)^{n+1}}$, where $A(t)$ is the unique integral
palindromic poly\-nomial of degree $c=k_X+n+1$ (the coindex) such that
$P_I(t)$ equals the series $P_X(t)$ up to and including degree
$\rd{\frac c2}$. If $c<0$ then $P_I=0$.

\item Each {\em orbifold part} for $Q\in\sB$ of type
$\frac1r(a_1,\dots,a_n)$ is of the form
$P\orb(Q,k_X)=\frac{B(t)}{(1-t)^n(1-t^r)}$, with
\begin{equation}
\label{eq!Bt}
B(t)=\InvMod\Bigl(\prod_{i=1}^n\frac{1-t^{a_i}}{1-t},
\frac{1-t^r}{1-t}, \rd{\frac c2}+1\Bigr)
\end{equation}
the unique Laurent polynomial supported in the interval
\begin{equation}
\Bigl[\ru{\tfrac{c-1}2}+1,\rd{\tfrac{c-1}2}+r-1\Bigr]
\label{eq!symI}
\end{equation}
equal to the inverse of\/ $\prod_{i=1}^n\frac{1-t^{a_i}}{1-t}$ modulo
$\frac{1-t^r}{1-t}$. The polynomial $B(t)$ has integral coefficients and
is palindromic of degree $k_X+n+r$.
\end{enumerate}
\end{theorem}
The quantities $\ru{\frac{c-1}2}=\rd{\frac c2}$ cause little headaches
of notation and computation. More conceptually, the numerator has
Gorenstein symmetry of degree $k+n+r$, giving the part as a whole
Gorenstein symmetry of degree $k$, and is a smallest residue modulo
$\frac{1-t^r}{1-t}$. The interval \eqref{eq!symI} as written is
manifestly symmetric, centred at $\frac{k+n+r}2$, and of length $\le
r-2$; it contains $r-1$ consecutive integers if $c$ is odd, and $r-2$ if
$c$ is even. The support of the numerator occupies the end points of the
stated interval in about half the cases. Compare Exercise~\ref{exc!ice}.

\begin{add} \label{ad}
We suppose $(X,D)$ is as in Theorem~\ref{th!main}, but relax the
projectively Gorenstein assumption to assume only that $K_X$ is
$\Q$-Cartier and numerically equivalent to $k_XD$. (In other words, omit
the projectively Cohen--Macaulay requirement.) Then the Hilbert series
of $X$ is
\begin{equation}
P_X(t)=J(t)+P_I(t)+\sum_{Q\in\sB} P\orb(Q,k_X)(t),
\end{equation}
with $P_I$ and $P\orb$ as above, where $J(t)=\sum j_mt^m$ is a
polynomial treating the irregularity of $\Oh_X(mD)$, with coefficients
\begin{equation}
\begin{aligned}
j_m &= h^0(\Oh_X(mD))+(-1)^nh^n(\Oh_X(mD))-\chi(\Oh_X(mD)) \\
 &= -\sum_{i=1}^{n-1} (-1)^i h^i(\Oh_X(mD)).
\end{aligned}
\end{equation}
In characteristic zero (or if some form of Kodaira vanishing holds) then
$J(t)$ has degree $\le k_X$.
\end{add}

\begin{exa} \rm
Consider the general hyper\-surface $X_{10}\subset\PP^4(1,1,2,2,3)$ with
coordinates $x_1,x_2,y_1,y_2,z$. Then $X_{10}$ is a 3-fold with
$5\times\frac12(1,1,1)$ orbifold points along $\PP^1_{\Span{y_1,y_2}}$
and a $\frac13(1,2,2)$ point at $P_z=(0,0,0,0,1)$. It has canonical
weight $k_X=1$ and coindex $c=k_X+n+1=5$. The Hilbert series is as
follows: the initial part
\begin{equation}
P_I=\frac{1-2t+3t^2+3t^3-2t^4+t^5}{(1-t)^4}
=1+t+\frac{t+t^2}{(1-t)^2}+2\frac{t^2+t^3}{(1-t)^4},
\end{equation}
takes care of $P_1=2$, $P_2=5$. The orbifold parts
\begin{equation}
P\orb(\textstyle\frac12(1,1,1),1) = \frac{-t^3}{(1-t)^3(1-t^2)}, \quad
P\orb(\textstyle\frac13(1,2,2),1) = \frac{-t^3-t^4}{(1-t)^3(1-t^3)}
\end{equation}
take care of the periodicity, giving
\begin{equation}
P_I+5\times P\orb(\tfrac12(1,1,1),1)+P\orb(\tfrac13(1,2,2),1)
=\frac{1-t^{10}}{(1-t)^2(1-t^2)^2(1-t^3)}.
\notag
\end{equation}

Here the numerator of $P_I$ is palindromic of degree $c=5$, so that
$P_I$ is Gorenstein symmetric of degree~1. The two $P\orb$ parts are
also integral and Gorenstein symmetric of degree~1, and they start with
$t^3$, so do not affect the first two plurigenera $P_1$ and $P_2$.
\end{exa}

\begin{caution} \label{ca!gl} \rm The initial part $P_I$ handles the
first plurigenera $P_1,\dots,P_{\rd{\frac c2}}$, but is not the {\em
leading term} of the Hilbert function controlling the order of growth of
the plurigenera: in this example $X_{10}\subset\PP^4(1,1,2,2,3)$ is a
canonical 3-fold with $K_X=\Oh_X(1)$, of degree
$K_X^3=\frac{10}{2\times2\times3} =\frac56$, whereas $P_I$ on its own
would correspond to degree $K^3=4$ (for this, sum the coefficients in
the numerator of $P_I$). In our formula, the orbifold parts contribute
to the global order of growth of the plurigenera, in this case
$5\times-\frac12$ and $-\frac23$.
\end{caution}

\subsection{Appendix: Symmetric integral polynomials}
\label{appendix!Hilb}
The shape of our Hilbert series in the nonsingular case comes directly
from the following result applied to Hilbert polynomials.

\begin{proposition} \label{p!binom}
\begin{enumerate} \renewcommand{\labelenumi}{(\Roman{enumi})}
\item Let $\sum_{m\ge0}\rho_m t^m$ be a power series, and assume that
$\rho_m=F(m)$ for all $m\ge m_0$, where $F(x)$ is a polynomial of degree
$n$ and $m_0\ge0$ an integer. Then $(1-t)^{n+1}\bigl(\sum_{m\ge0} \rho_m
t^m\bigr)$ is a polynomial in $t$ of degree $\le m_0+n+1$.

\smallskip
\item Let $F(x)\in\Q[x]$ be a polynomial taking integer values $F(m)$
for all $m\in\Z$. Then $F$ is an integral linear combination of the
binomial coefficients:
\begin{equation}
F(x) = \sum_{\nu=0}^n c_\nu \binom x\nu \quad\hbox{with $c_\nu\in\Z$}.
\end{equation}
Here $n=\deg F$, and there are $n+1$ integral coefficients $c_\nu$ to
specify.

\smallskip
\item Let $F(x)\in\Q[x]$ be a polynomial taking integer values $F(m)$
for all $m\in\Z$. Assume that $F$ satisfies $(-1)^n F(k-x)\equiv F(x)$
for an integer $k$, where $\deg F=n$. Then $F(X)$ and its associated
power series $\sum_{m\ge0}F(m)t^m$ are integral linear combinations of
standard terms as follows
\begin{description}
\smallskip
\item[if $n\equiv k+1\mod 2$]
\begin{equation}
\begin{aligned}
F(x) &= \sum_{\substack{-k\le\nu\le n \\ \nu \equiv n\mod 2}}
b_\nu \binom{x+\frac{\nu-k-1}2}{\nu}, \\
\sum_{m\geq 0}F(m)t^m &=
\sum_{\substack{-k\le\nu\le n \\ \nu \equiv n\mod 2}} b_\nu
\frac {t^{\frac{\nu+k+1}2}}{(1-t)^{\nu+1}}\,.
\end{aligned}
\end{equation}

\smallskip
\item[or if $n\equiv k\mod 2$]
\begin{equation}
\begin{aligned}
F(x) &= \sum_{\substack{-k\le\nu\le n \\ \nu \equiv n\mod 2}}
b_\nu \left( \binom{x+\frac{\nu-k}2}{\nu} + \binom{x+\frac{\nu-k-2}2}{\nu}\right),
\\
\sum_{m\geq 0}F(m)t^m &=
\sum_{\substack{-k\le\nu\le n \\ \nu \equiv n\mod 2}} b_\nu
\frac {(1+t)t^{\frac{\nu+k}2}}{(1-t)^{\nu+1}}\,.
\end{aligned}
\end{equation}
\end{description}
There are $\rd{\frac {k+n+1}2}$ integral coefficients $b_\nu$ to specify.
\end{enumerate}
\end{proposition}

In part (II) or (III), it is enough to assume that $F(m)\in\Z$ or
$F(m)\in\Z$ and $(-1)^{n}F(k-m)=F(m)$ for all $m$ in an interval of
length $n+1$. The proof is a little exercise. Hint: Use induction based
on $F(x)-F(x-1)$. \qed

\section{Ice cream functions} \label{s!ice}

\subsection{Fun calculation} \label{ss!ice}
``Income $\frac37$ per day means ice cream on Wednesdays, Fridays and
Sundays.'' Consider the step function $i\mapsto \rd{\frac{3i}7}$, where
$\rd{\,}$ denotes the rounddown or integral part. As a Hilbert series, it
gives
\begin{equation}
\label{eq!ice37}
P(t)=\sum_{i\ge0} \rd{\frac{3i}7}
t^i=0+0t+0t^2+t^3+t^4+2t^5+2t^6+3t^7+\cdots,
\end{equation}
with closed form
\begin{equation}
P(t)=\frac{t^3+t^5+t^7}{(1-t)(1-t^7)}.
\label{eq!357}
\end{equation}
Indeed, $\rd{\frac{3i}7}$ increments by~1 when $i=0,3,5$ modulo 7, so
that
\begin{equation}
\label{eq!icerow}
(1-t)P(t)=t^3+t^5+t^7+t^{10}+\cdots
\end{equation}
is the sum over the jumps, that repeat weekly. Multiplying
\eqref{eq!icerow} by $1-t^7$ cuts the series down to the first week's
ice cream ration:
\begin{equation}
(1-t)(1-t^7)P(t)=t^3+t^5+t^7.
\end{equation}

The numerator $t^3+t^5+t^7$ can be seen as the
\begin{equation}
\hbox{inverse of }\frac{1-t^5}{1-t}=1+t+t^2+t^3+t^4 \mod
\frac{1-t^7}{1-t}=1+t+t^2+t^3+t^4+t^5+t^6.
\notag
\end{equation}
Indeed, long multiplication gives
\begin{multline}
\label{eq!lmult}
(1+t+t^2+t^3+t^4)\times(t^3+t^5+t^7)= \\
\setcounter{MaxMatrixCols}{30}
\renewcommand{\arraycolsep}{.075em}
\begin{matrix}
&&&&&&&t^3&+&t^4&+&t^5&+&t^6&+&t^7&+& \\
&&&&&&&&&&&t^5&+&t^6&+&t^7&+&t^8&+&t^9 \\
&&&&&&&&&&&&&&&t^7&+&t^8&+&t^9&+&t^{10}&+&t^{11}\\[4pt]
&&&&&&=&t^3&+&t^4&+&2t^5&+&2t^6&+&3t^7&+&2t^8&+&2t^9&+&t^{10}&+&t^{11}\\[4pt]
\equiv&3&+&2t&+&2t^2&+&2t^3&+&2t^4&+&2t^5&+&2t^6 &
\equiv&1,
\end{matrix}
\end{multline}
where $\equiv$ denotes congruence modulo $\frac{1-t^7}{1-t}$. Here
$5=\InvMod(3,7)$ is the inverse of 3 modulo 7. The product in
\eqref{eq!lmult} has $5\times3=15\equiv1\mod 7$ terms that distribute
themselves equitably among the 7 congruence classes, except that $t^7$
appears once for each of the 3 terms in the second factor.

There are several other meaningful expressions for $P(t)$. Under
$\equiv$, the bounty $t^3+t^5+t^7$ can be viewed as famine
$-t-t^2-t^4-t^6$ ``no ice cream on Mondays, Tuesdays, Thursdays or
Saturdays''. In other words,
\begin{equation}
P(t)=\frac{t^3+t^5+t^7}{(1-t)(1-t^7)}
=\frac t{(1-t)^2}+\frac{-t-t^2-t^4-t^6}{(1-t)(1-t^7)}.
\label{eq!Pt}
\end{equation}
Because $t^7\equiv1$, we can shift the exponents of $t$ up or down by 7:
\begin{equation}
\frac{t^{-4}+t^{-2}+1}{(1-t)(1-t^7)}\quad\hbox{or}\quad
\frac{-t^{-1}-t-t^2-t^4}{(1-t)(1-t^7)}
\end{equation}
so ``ice cream rations from Monday before the start of term'' or
``famine from the previous Saturday''. More generally, working mod
$\frac{1-t^r}{1-t}$, we can shift exponents $t^b\mapsto t^{b-ir}$ mod
$r$ and subtract a multiple of $1+t+\cdots+t^{r-1}$ to {\em fold} any
Laurent polynomial into any desired interval $[t^a,\dots,t^{a+r-2}]$ of
length $r-2$. Of these possible shifts as Laurent polynomials with
short support, $t^{7i}(t^3+t^5+t^7)$ is palindromic of degree $10+14i$,
and $t^{7i}(-t^{-1}-t-t^2-t^4)$ is palindromic of degree $3+14i$, {\em
and no other}.

In ``macroeconomic'' terms, the order of growth is the linear function
$\frac{3i}7$ with fractional seasonal corrections, that is,
\begin{equation}
P(t)=\frac37\cdot\frac t{(1-t)^2}
+\frac{-\frac37t-\frac67t^2-\frac27t^3-
\frac57t^4-\frac17t^5-\frac47t^6}{1-t^7}
\label{eq!Dkd0}
\end{equation}
(``on Mondays, we lose $\frac37$ in small change'', etc.). Notice the
coefficient $\frac17$ of $t^5$: the inverse of 3 modulo 7 is 5, so as we
enjoy our second ice cream of the week on Fridays, we lose $\frac17$,
the unit of small change.

We can average out the seasonal corrections in \eqref{eq!Dkd0} to sum to
zero, giving
\begin{equation}
P(t)=\frac37\cdot\frac{2t-1}{(1-t)^2}
+\frac{\frac37-\frac37t^2+\frac17t^3-\frac27t^4+\frac27t^5-\frac17t^6}{1-t^7},
\label{eq!Dkd}
\end{equation}
where the coefficients $\frac37,0,-\frac37,\frac17,\frac27,\frac27,
-\frac17$ are the Dedekind sums $\si_i\bigl(\frac17(5)\bigr)$ (see
Definition~\ref{df!Dkd} and compare \cite{YPG}, Theorem~8.5).

The expressions \eqref{eq!357}--\eqref{eq!Dkd} represent different views
on numerical functions that grow with periodic corrections. Our main aim
is to explain the orbifold contributions $P\orb$ in
Theorem~\ref{th!main} as minor variations on this simple-minded material.

\begin{exc}[One dimensional ice cream functions] \label{exc!ice} \rm
Let $0<a<r$ be coprime integers, $k\equiv-a$ mod $r$ and set
$b=\InvMod(a,r)$. An integer is an {\em ice cream day} for $\frac1r(a)$
if its congruence class is one of the $b$ distinct classes
\begin{equation}
\{0,a,2a,\dots,a(b-1)\} \mod r,
\end{equation}
and a {\em non ice cream day} if its congruence class is one of the
$r-b$ classes
\begin{equation}
\{1,a+1,2a+1,\dots,a(r-b-1)+1\} \mod r.
\end{equation}
The two sets are complimentary because $ai+1\equiv a(i+b)$ mod $r$.

Then the numerator $B(t)$ of \eqref{eq!Bt} is one of
\begin{equation}
(\I) = \sum  t^j \quad
 \hbox{summed over ice cream days in the interval \eqref{eq!symI}}
\label{eq!I}
\end{equation}
or
\begin{equation}
(\II) = \sum -t^j \quad
 \hbox{summed over non ice cream days in \eqref{eq!symI}}.
\label{eq!II}
\end{equation}
\paragraph{\bf Hint} Since $ab\equiv1$ mod $r$, the inverse of
$\frac{1-t^a}{1-t}$ mod $1-t^r$ is
$\frac{1-t^{ab}}{1-t^a}=\sum_{i=0}^{b-1} t^{ai}$.
\end{exc}

\subsection{The function Inverse Mod} We start with the following basic
result.

\begin{theorem} \label{th!invm}
Fix an integer $\ga$ and a monic polynomial $F\in\Q[t]$ of degree $d$
with nonzero constant term.
\begin{enumerate}
\renewcommand\labelenumi{(\Roman{enumi})}
\item The quotient ring $\Q[t]/(F)$ is a $d$-dimensional vector space
over $\Q$ and $t$ is invertible in it, so that
$\Q[t]/(F)=\Q[t,t^{-1}]/(F)$.

\item Any range $[t^\ga,\dots,t^{\ga+d-1}]$ of $d$ consecutive Laurent
monomials maps to a $\Q$-basis of $\Q[t]/(F)$.

\item If $A\in\Q[t]$ is coprime to $F$, we can write its inverse modulo
$F$ uniquely as a Laurent polynomial $B$ with support in
$[t^\ga,\dots,t^{\ga+d-1}]$.

\end{enumerate}
\end{theorem}
\paragraph{\bf Proof}
This is all trivial. The leading term of $F$ is nonzero, so
$1,t,\dots,t^{d-1}$ base $\Q[t]/(F)$. The constant term of $F$ is
nonzero so $t$ is coprime to $F$, hence invertible modulo $F$.
Multiplication by $t$ is an invertible linear map, so multiplication by
$t^\ga$ for any $\ga\in\Z$ takes a basis to another basis. If $A$ is
coprime to $F$ it is invertible in $\Q[t]/(F)$, and its inverse has a
unique expression in any basis. \qed\par\medskip

\begin{dfn} \rm For coprime polynomials $A,F\in \Q[t]$ we set
\begin{equation}
\InvMod(A,F,\ga)=B
\end{equation}
with $B$ as in (III). That is, $B\in\Q[t,t^{-1}]$ is the uniquely
determined Laurent polynomial supported in $[t^\ga,\dots,t^{\ga+d-1}]$
with $AB\equiv1\mod F$. For different $\ga\in\Z$, these inverses are
congruent modulo $F$, but of course different polynomials in general. We
also write $\InvMod(A,F)$ with unspecified support for any inverse of
$A$ modulo $F$ in $\Q[t]$.
\end{dfn}

Fix positive integers $r$ and $a_1,\dots,a_n$ and set
\begin{equation}
\label{eq!AF}
A=\prod_{j=1}^n(1-t^{a_j}), \quad h=\hcf(1-t^r,A) \quad\hbox{and}\quad
F = \frac{1-t^r}{h}.
\end{equation}
The polynomial $F$ is the monic polynomial with simple roots only at the
$r$th roots of unity with $A(\ep)\ne0$, or equivalently $\ep^{a_j}\ne1$
for all $a_j$. Since we take out the hcf, $A$ and $F$ are coprime.
Theorem~\ref{th!invm} applies to give $\InvMod(A,F,\ga)$, the inverse of
$A$ modulo $F$ with support in $[t^\ga,\dots,t^{\ga+d-1}]$, where
$d=\deg F$ and $\ga\in\Z$ is arbitrary.

We show how to compute $\InvMod$:
\begin{alg} \label{algInvMod} \rm
If $\ga\ge0$ then $t^\ga A$ and $F$ are coprime polynomials. Set $d=\deg
F$. The Euclidean algorithm in $\Q[t]$ provides a unique solution to
\begin{equation}
t^\ga AB+FG=1,
\label{eq!Eucl}
\end{equation}
with $B\in\Q[t]$ a polynomial of degree $<d$. Then
$\InvMod(A,F,\ga)=t^\ga B$.

If $\ga<0$, choose some $m$ with $mr+\ga\ge0$, and solve
\begin{equation}
t^{mr+\ga}AB+FG=1
\end{equation}
by the Euclidean algorithm. Then $\InvMod(A,F,\ga)
=t^{mr+\ga}B/t^{mr}=t^\ga B$. This trick works because
$t^{mr}\equiv1\mod F$. For more general poly\-nomials $F$, one would
need to calculate powers of the matrix $M_t$ corresponding to
multiplication by $t$ in $\Q[t]/(F)$; in our case, $M_t^r=1$.
\end{alg}

The isolated case is when $a_1,\dots,a_n$ are coprime to $r$, so $h=1-t$
and $F=1+t+\cdots+t^{r-1}$ has degree $d=r-1$ and roots
$\ep\in\bmu_r\setminus\{1\}$. If moreover $r$ is prime then $F$ is the
cyclotomic polynomial, and working modulo $F$ is essentially the same
thing as setting $t=\ep$ a primitive $r$th root of unity.

\subsection{Dedekind sums as Inverse Mod} \label{s!DIM}

We now recall Dedekind sums, and relate them to the function $\InvMod$.
\begin{dfn} \label{df!Dkd} \rm
We define the $i$th {\em Dedekind sum} $\si_i$ by
\begin{equation}
\si_i\left(\tfrac1r(a_1,\dots,a_n)\right)
=\frac1r\sum_{\substack{\ep\in\bmu_r \\ \ep^{a_j}\ne1\,\forall j = 1,\dots,n}}
\frac{\ep^i}{(1-\ep^{a_1})\cdots(1-\ep^{a_n})},
\label{eq!defsi}
\end{equation}
where $\ep$ runs over the $r$th roots of unity for which the denominator
is nonzero. Proposition~\ref{p!Dkd} characterizes the $\si_i$ as
solutions to a set of linear equations. We combine them into the {\em
Dedekind sum polynomial}
\begin{equation}
\De\bigl(\tfrac1r(a_1,\dots,a_n),t\bigr)=\sum_{i=1}^r\si_{r-i}t^i
\quad\hbox{with support in $[t,\dots,t^r]$}.
\label{eq!defDe}
\end{equation}

It is obvious that $\si_i=\si_{r+i}$. Therefore we only need to consider
$\si_i$ for $i=0,1,\dots,r-1$. In the coprime case, the sum in
\eqref{eq!defsi} runs over all nontrivial $r$th roots of unity. To
stress that $a_1,\dots,a_n$ are not all coprime to $r$, we may call
$\si_i$ the $i$th {\em generalized} Dedekind sum.
\end{dfn}

\begin{proposition} \label{p!Dkd}
Consider the $r\times r$ system of linear equations
\begin{equation}
\sum_{i=0}^{r-1} \si_i \ep^i =
\begin{cases}
0 & \hbox{if $\ep\in \bmu_{a_j}$ for some $j$,} \\
\dfrac1{(1-\ep^{-a_1})\cdots(1-\ep^{-a_n})} & \hbox{otherwise.}
\end{cases}
\label{eq!Dk1}
\end{equation}
in unknowns $\si_i$ indexed by $i\in\Z/r=\Hom(\bmu_r,\C^\times)$, with
equations indexed by $\ep\in\bmu_r$.

Then \eqref{eq!Dk1} is a nondegenerate system, with unique solution
the Dedekind sums $\si_i=\si_i\left(\tfrac1r(a_1,\dots,a_n)\right)$.
\end{proposition}

\paragraph{\bf Proof} Fix a prime root of unity $\ep\in\bmu_r$. Then
$(\ep^{ij})_{i,j=0,\dots,r-1}$ is a Vandermonde matrix, with inverse
$\tfrac1r \left(\ep^{-ij}\right)_{i,j=0,\dots,r-1}$. \qed\par\medskip

\begin{lemma} \label{l!coset}
The polynomial $\De$ \eqref{eq!defDe} is divisible by
$h=\hcf(A,1-t^r)$.
\end{lemma}

\paragraph{\bf Proof} The roots of $h$ are the $\ep\in\bmu_r$ for which
$\ep^{a_j}=1$ for some $j$. An equivalent statement is that if $\be$ is
a common divisor of $r$ and some $a_j$ then
\begin{equation}
\sum_{\substack{i=0,\dots,r-1,\\ i\equiv d\mod \be}} \si_i=0
\quad\hbox{for any integer $d$.}
\end{equation}
In words, the average of the $\si_i$ over any coset of
$\be\Z/r\subset\Z/r$ is zero. In particular, $\sum_{i=0}^{r-1} \si_i=0$.

Note that $\ep\in\bmu_r$ gives $\ep^{\be}\in\bmu_{r/\be}$. Then by
Definition~\ref{df!Dkd},
\begin{multline}
\si_d+\si_{d+\be}+\cdots+\si_{d+r-\be}= \\
\frac1r\sum_{\substack{\ep\in\bmu_r \\ \ep^{a_j}\ne1\,\forall j}}
\frac{\ep^d}{\prod_j(1-\ep^{a_j})}\left(1+\ep^{\be}+\ep^{2\be}+\cdots+
\ep^{\be\left(\frac r\be-1\right)}\right)=0.\notag \qed
\end{multline}
For example
\begin{equation}
\si_i\bigl(\tfrac1{14}(1,2,5,7)\bigr)
=\frac1{14}\bigl\{-2,-2,-1,\tfrac12,0,-\tfrac12,1,2,2,1,-\tfrac12,0,\tfrac12,-1\bigr\},
\end{equation}
with $\si_i+\si_{7+i}=\sum_{l=0}^6 \si_{2l+i}=\sum_{l=0}^{13}
\si_{l+i}=0$ for each $i$.

The next result was first stated and proved by Buckley \cite[Theorem
2.2]{eurosim}, following the ideas of \cite{YPG}.

\begin{theorem} \label{th!eurosim}
Let $A$, $h$ and $F$ be as in \eqref{eq!AF} and $\De$ as in \eqref{eq!defDe}.
Then
\begin{equation}
\De = ht \cdot \InvMod(htA,F,0).
\label{eq!Si}
\end{equation}
\end{theorem}

\medskip
\paragraph{\bf Proof}

Since $F$ is coprime to $t$ and $h$, multiplying by $th$ before and
after taking InverseMod does not change the result modulo~$F$. The factor
$h$ makes the right-hand side of \ref{eq!Si} divisible by $h$ in
accordance with Lemma~\ref{l!coset}. The factor $t$ then folds it from a
polynomial supported in $[0,\dots,r-1]$ to $[1,\dots,r]$, as in
\eqref{eq!Eucl}. It only remains to prove that
\begin{equation}
\De\equiv\InvMod(A,F,\ga) \in\Q[t]/(F)
\quad\hbox{for any $\ga$},
\end{equation}
or equivalently, that $B(t):=A(t)\De$ is congruent to $1$ modulo $F$.
For this, substitute any root $t=\ep$ of $F$ in $B$ and use
\eqref{eq!Dk1} with the inverse value of $\ep$. This gives
\begin{equation}
B(\ep) = A(\ep)\De = \frac{A(\ep)}{\prod_j(1-\ep^{a_j})} = 1.
\end{equation}
This holds for every root $\ep$ of $F$, so $B(t)-1$ is divisible by $F$,
that is,
\begin{equation}
A(t)\De\equiv1 \mod F. \qed
\end{equation}
\par\medskip

\begin{proposition} \label{p!np}
Assume all the $a_i$ are coprime to $r$, so that $h=1-t$,
$F=\frac{1-t^r}{1-t}$ and $d=\deg F=r-1$. Then for any $\ga$,
\begin{align}
(1-t)^n\De &\equiv
\InvMod\Bigl(\prod_{j=1}^n\frac{1-t^{a_j}}{1-t},F\Bigr)
\label{eq!prod} \\
& \equiv \InvMod\Bigl(\frac A{(1-t)^n},F,\ga+1\Bigr)
=\sum_{l=\ga+1}^{\ga+r-1} \theta_lt^l,
\notag
\end{align}
with integer coefficients $\theta_l=\sum_{s=0}^n (-1)^s \binom ns
(\si_{s-l}-\si_{s-\ga})\in\Z$
\end{proposition}

\medskip
\paragraph{\bf Proof}
Replace the InvMod of a product by the product of InvMods. Each factor
$\InvMod\bigl(\frac{1-t^{a_j}}{1-t},F,1\bigr)$ is a polynomial with
integral coefficients; indeed, by the calculation of \ref{ss!ice}, or
Exercise~\ref{exc!ice}, it is the ice cream function for $\frac{b_j}r$
where $b_j=\InvMod(a_j,r)$.
\qed\par

\begin{exc}[Serre duality, Gorenstein symmetry]\label{exc!ser} \rm
If $X$ is projectively Gorenstein of canonical weight $k_X$, prove the
following:

\begin{enumerate}

\item Each $Q=\frac1r(a_1,\dots,a_n)$ satisfies $k_X+\sum_{j=1}^n
a_j\equiv0\mod r$.

\item The $\si_i$ are $(-1)^n$ symmetric under $i\mapsto \sum a_j-i$.
[Hint: replace $\ep\mapsto\ep^{-1}$ in the characterization
\eqref{eq!Dk1} of the $\si_i$, or in \eqref{eq!defsi}.]

\item Now let $\theta_l$ be as in Proposition~\ref{p!np}. Then
$l_1+l_2\equiv k_X+n\mod r$ implies $\theta_{l_1}=\theta_{l_2}$. In
particular, for $c$ even and $\ga=\frac c2$, we have
$\theta_{\ga+r-1}=0$, since $\theta_{\ga}=0$ by definition.
\end{enumerate}
\end{exc}

\subsection{Ice cream gives the correct periodicity} \label{ss!Porb}
There are two expressions for the orbifold contributions to RR. The
first, given in \cite{YPG}, is in terms of Dedekind sums:
\begin{equation}
\frac{\sum_{i=1}^{r-1}(\si_{r-i}-\si_0)t^i}{1-t^r}.
\end{equation}
The alternative introduced here is the {\em ice cream function}
\begin{equation} \label{PorbInteg!eq}
P\orb\bigl(\tfrac1r(a_1,\dots,a_n),k_X\bigr)=
\frac{B(t)}{(1-t)^n(1-t^r)},
\end{equation}
with
\begin{equation}
B(t)=\InvMod\Bigl(\prod_{i=1}^n\frac{1-t^{a_i}}{1-t},
\frac{1-t^r}{1-t}, \rd{\frac c2}+1\Bigr)
\end{equation}
as in \eqref{eq!Bt}. The first is strictly periodic (because of the
denominator $1-t^r$), but fractional. The second is integral by
Proposition~\ref{p!np}, and Gorenstein symmetric of degree $k$, but has
order of growth $O(m^n)$. They both give the same periodicity, as a
simple consequence of Proposition~\ref{p!np}. The point already appeared
clearly in the different treatments of $P(t)$ in \eqref{eq!ice37} and
\eqref{eq!Dkd0}--\eqref{eq!Dkd}.

\begin{cor} \label{cor!per}
\begin{equation}
P\orb(\tfrac1r(a_1,\dots,a_r),k_X)-\frac{\sum_{i=1}^{r-1}(\si_{r-i}-\si_0)t^i}{1-t^r}
=\frac{C(t)}{(1-t)^{n+1}}
\end{equation}
with numerator $C(t)\in\Q[t]$.
\end{cor}
Indeed, put the left hand side over the common denominator $(1-t)^n(1-t^r)$ and
use Proposition~\ref{p!np}. \qed\par\medskip

In the noncoprime case, the result is similar, and we leave the proof as
an exercise.

\begin{proposition} \label{p!np1}
Set $s_i=\hcf(a_i,r)$, so that $\hcf(1-t^{a_i},1-t^r)=1-t^{s_i}$, and
write $d=\deg F$. Then for any $\ga$,
\begin{align}
\prod(1-t^{s_i})\cdot\De &\equiv
\InvMod\Bigl(\prod_{j=1}^n\frac{1-t^{a_j}}{1-t^{h_j}},F,\ga+1\Bigr)
\label{eq!prod1} \\
& = \InvMod\Bigl(\frac A{\prod(1-t^{s_i})},F,\ga+1\Bigr)
=\sum_{l=\ga+1}^{\ga+d} \theta_lt^l,
\notag
\end{align}
with integer coefficients $\theta_l$.

\rm For appropriate choice of $\ga+1$, this gives a {\em generalized ice
cream function} that is integral and Gorenstein symmetric, having the
same $r$ periodicity as $\frac{\De}{1-t^r}$. Compare \eqref{eq!genIce}.
\end{proposition}

\begin{exa} \label{exa!57} \rm
The ice cream function of \ref{ss!ice} corresponds to
$\si_i\left(\frac17(5)\right)$: the periodic rounding loss of
\eqref{eq!Dkd0}--\eqref{eq!Dkd} is
\begin{align*}
\sum_{i=0}^6 \si_{7-i}t^i &= \tfrac17(3-3t^2+t^3-2t^4+2t^5-t^6) \\
& \equiv \InvMod\bigl(1-t^5,\tfrac{1-t^7}{1-t},0\bigr).
\end{align*}
Multiplication by $1-t$ gives a Gorenstein symmetric polynomial with
integral coefficients $\theta_l$
\begin{align*}
(1-t)&\times\tfrac17(3-3t^2+t^3-2t^4+2t^5-t^6) \\
& \qquad \equiv t^3+t^5+t^7
=\InvMod\Bigl(\frac{1-t^5}{1-t},\frac{1-t^7}{1-t},3\Bigr).
\end{align*}

The fractional divisor $\frac37P$ on a nonsingular curve is an orbifold
point of type $\frac17(5)$, with orbinate in $\Oh(5)$ having genuine pole
of order two, but fractional zero of order $\frac17$ in ``lost change''.
As we saw in Exercise~\ref{exc!ice}, the same considerations apply with
$\frac37$ replaced by a general reduced fraction $\frac ar$,
corresponding to the orbifold point $\frac1r(b)$ with $b=\InvMod(a,r)$.

Consider for example the weighted projective line $X=\PP(5,7)$. It has
$k_X=-12$, and has two orbifold points of type $\frac17(5)$ and
$\frac15(2)$. Its Hilbert series
\begin{equation}
P_X(t)=\frac1{(1-t^5)(1-t^7)}
\end{equation}
satisfies Theorem~\ref{th!main}: since $c=-10<0$, the initial part
$P_I=0$. Then
\begin{align*}
P_X(t)&=P\orb \bigl(\tfrac17(5),-12\bigr)+P\orb \bigl(\tfrac15(2),-12\bigr) \\
&=\frac{t^{-4}+t^{-2}+1}{(1-t)(1-t^7)}+\frac{-t^{-4}-t^{-2}}{(1-t)(1-t^5)},
\end{align*}
where $-t^{-4}-t^{-2}=\InvMod(\frac{1-t^2}{1-t},\frac{1-t^5}{1-t},-4)$.
\end{exa}

\begin{exc} \label{exc!57} \rm
Fun and games with the ice cream functions of \ref{ss!ice}.
\begin{enumerate}
\item An elliptic curve polarized by $A=\frac37P$ embeds as
$C_{15}\subset\PP(1,5,7)$ with canonical weight 2, that is,
$K_{C,\mathrm{orb}}=2A=\frac67P$.

\item A quasismooth complete intersection $C_{10,15}\subset\PP(1,3,5,7)$
is a curve of genus 7 with $K_C=3P+9Q$ having $P$ as an orbifold point
of type $\frac17(5)$, polarized by $A=\frac37P+Q$ and having
$K_{C,\mathrm{orb}}=9A$. (Its initial part $P_I$ is quite involved.)

\item A curve of genus 2 polarized by $P+\frac37Q$ with $P$ a Weierstrass
point embeds in $\PP(1,2,3,5,7)$ as a Pfaffian with Hilbert numerator
\begin{equation}
1-t^6-t^7-t^8-t^9-t^{10}+t^{10}+t^{11}+t^{12}+t^{13}+t^{14}-t^{20}.
\notag
\end{equation}

\end{enumerate}
\end{exc}

\section{Proof of Main Theorem}
\label{s!pf}

\subsection{The existence of the Riemann--Roch formula} \label{s!RRD}

Let $X$ be a normal projective $n$-fold; assume that the singularities
of $X$ are isolated, rational and $\Q$-factorial. We want to calculate
$\chi(\Oh_X(D))$ for $D$ a Weil divisor on $X$ using the RR formula
\begin{equation}
\label{eq!RR0}
\bigl(\ch(\Oh_X(D))\cdot\Td(T_X)\bigr)[n],
\end{equation}
that is, the component of top degree $n$ of the product of
\begin{align*}
\ch(\Oh_X(D))&=\exp(D)=\sum \tfrac{D^i}{i!}, \quad \hbox{and} \\
\Td(T_X)&= \sum_{i=0}^n \Td_i(T_X) \\
&=1-\tfrac12K_X+\tfrac1{12}(K_X^2+c_2)-\tfrac1{24}K_X c_2 \\
& \qquad -\tfrac1{720}(K_X^4-4K_X^2 c_2-3 c_2^2+K_X c_3+c_4)+\cdots.
\end{align*}
We must get around the problem that the terms in \eqref{eq!RR0} are not
defined, because $T_X$ is not a vector bundle on a singular $X$. For
this, we use the following conventions. First, choose a resolution of
singularities $f\colon Y\to X$ that is an isomorphism over the
nonsingular locus of $X$.
\begin{enumerate}
\renewcommand{\labelenumi}{(\alph{enumi})}
\item Replace the degree $n$ term $\Td_n(T_X)$ in \eqref{eq!RR0}
by $\chi(\Oh_X)=\chi(\Oh_Y)=\Td_n(T_Y)$.
\item Replace the terms involving a product with $D$ on $X$ by the
same expression on $Y$ involving its pullback as a $\Q$-Cartier
divisor. In more detail: the pullback of a $\Q$-Cartier divisor $D$ is
defined as usual by $f^*D=\frac1mf^*(mD)$ with $mD$ Cartier. Except for
$\Td_n(T_X)$, the terms in \eqref{eq!RR0} are $D^i\Td_{n-i}(T_X)/i!$
with $i\ge1$, and we replace
\[
D^i\Td_{n-i}(T_X) \quad\hbox{by}\quad
(f^*D)^i\Td_{n-i}(T_Y).
\]
\end{enumerate}

\begin{remark} \rm
Our interpretation of \eqref{eq!RR0} is independent of the choice of the
resolution $Y$. Indeed, $\chi(\Oh_Y)$ is a birational invariant. Each of
the other terms involves a product with the $\Q$-Cartier divisor $D$;
now a multiple $mD$ is linearly equivalent to a linear combination of
nonsingular prime divisors disjoint from the singularities of $X$, so we
can calculate $D^i\Td(T_X)$ for $i>0$ on the nonsingular locus of $X$
itself.

In (a), we use $\chi(\Oh_X)$ as a substitute for $\Td_n(T_X)$. In the
3-fold case, it is well known that the expression $\Td_3(T_X)=-\frac1{24}
K_X\cdot c_2(T_X)$ can be defined using the same trick as in (b) (taking
the pullback of the $\Q$-Cartier divisor $K_X$), but {\em is not equal}
to $\chi(\Oh_X)$ in general. See \cite[Corollary~10.3]{YPG} and compare
Kawamata \cite[2.2]{Ka1} and \cite{Ka2}.
\end{remark}

\begin{theorem} \label{th!RRD}
Let $X$ be a normal projective $n$-fold with isolated, rational,
$\Q$-factorial singularities and $f\colon Y\to X$ as above. Then the
expression
\begin{align}
\label{eq!RRD}
\RR(D) &= \chi(\Oh_X) + \sum_{i=1}^n \frac1{i!}(f^*D)^i\Td_{n-i}(T_Y) \\
& = \hbox{``\,$\bigl( \ch(\Oh_X(D))\cdot \Td(T_X) \bigr)[n]$''} \notag
\end{align}
is a polynomial in the $\Q$-Cartier Weil divisor $D$ such that for every
$D$, the difference
\begin{equation}
\chi(X,\Oh_X(D)) - \RR(D) = \sum_{Q\in\Sing X} c_Q(D)
\end{equation}
is a sum of fractional terms $c_Q(D)\in\Q$ depending only on the local
analytic type of $X$ and $D$ at each singular point $Q$ of $X$.
\end{theorem}

\paragraph{\bf Plan of proof} We set $\sL=f^*\Oh_X(D)/\mathrm{torsion}$,
which is a torsion free sheaf of rank~1 on $Y$, and write
$\Oh_Y(D_Y)=\sL^{\vee\vee}$ for its reflexive hull, which is an
invertible sheaf. The proof has two parts: the first uses the Leray
spectral sequence to compare $\chi(X,\Oh_X(D))$ with
$\chi(Y,\Oh_Y(D_Y))$, given by RR on $Y$. After this, we compare the RR
formula for $D_Y$ on $Y$ with our interpretation $\RR(D)$ of the RR
formula for $D$ on $X$. No sooner said than done.

The reflexive hull of $\sL$ fits in a short exact sequence
\begin{equation}
\label{eq!DY}
0 \to \sL \to \Oh_Y(D_Y) \to \sQ \to 0,
\end{equation}
where the cokernel $\sQ$ has support of codimension $\ge2$ in $Y$ contained
in the exceptional locus of $f$.

Now $f_*\sL=f_*\Oh_Y(D_Y)=\Oh_X(D)$ because $\Oh_X(D)$ is saturated.
Moreover, all the sheaves $R^if_*\sL$ for $i\ge1$ and $R^if_*\sQ$ for
$i\ge0$ are finite dimensional vector spaces supported at the singular
points of $X$.

Now the Leray spectral sequence together with the long exact sequence
associated with \eqref{eq!DY} gives
\begin{align}
\chi(Y,\Oh_Y(D_Y)) &= \chi(\sL) + \chi(\sQ) \notag \\[-6pt]
&= \chi(X,\Oh_X(D)) + \sum_{i=1}^{n-1} (-1)^i h^0(X,R^if_*\sL) \\[-6pt]
&\qquad + \sum_{i=0}^{n-1} (-1)^i h^0(X,R^if_*\sQ). \notag
\end{align}
We deduce that $\chi(X,\Oh_X(D)) = \chi(Y,\Oh_Y(D_Y))+\sP$, where
\begin{equation} \label{eq!st1}
\sP = -\sum\nolimits_i(-1)^ih^0(X,R^if_*\sL)
-\sum\nolimits_i(-1)^ih^0(X,R^if_*\sQ)
\end{equation}

The second part of the proof depends on the exceptional locus of $f$.
Write $E_j$ for the exceptional divisors over the singular points, and
set
\begin{equation}
\label{eq!discD}
f^*D=D_Y+F, \quad\hbox{where}\quad
F=\sum\nolimits_jm_jE_j \quad\hbox{with $m_j\in\Q$.}
\end{equation}
The exceptional divisor $F$ here is the fixed part of the
birational transform of the linear system $|D+H|$ for any
sufficiently ample Cartier divisor $H$ on $X$.

Then
\begin{equation}
\label{eq!F}
\chi(Y,\Oh_Y(D_Y))-\RR(D) =
\sum_{i>0} \tfrac1{i!} (-F)^i\Td_{n-i}(T_Y).
\end{equation}

In fact, our interpretation $\RR(D)$ of $\ch(D)\cdot\Td(T_X)$ replaces
$\Td_n(T_X)$ by $\chi(\Oh_X)=\chi(\Oh_Y)=\Td_n(T_Y)$, and
$D^i\Td_{n-i}(X)$ by $(f^*D)^i\Td_{n-i}(Y)$, whereas the terms in RR on
$Y$ are $D_Y^i\Td_{n-i}(Y)$. Therefore the difference in \eqref{eq!F} is
\begin{equation}
\sum_{i>0} \tfrac1{i!} (D_Y^i - (f^*D)^i)\Td_{n-i}(T_Y).
\end{equation}
However, $f^*D$ is orthogonal to the exceptional divisors, and one
checks using the binomial expansion that
$D_Y^i-(f^*D)^i=(D_Y-f^*D)^i=(-F)^i$.

In conclusion, the difference required in Theorem~\ref{th!RRD} is
\begin{equation}
\chi(X,\Oh_X(D))-\RR(D)=\sP +
\sum_{i>0} \tfrac1{i!} (-F)^i\Td_{n-i}(T_Y).
\end{equation}

We can choose the resolution of singularities of $X$ and $D$ depending
only on the local analytic type of $X$ and $D$. The resolution
determines the sheaves $\sL$ and $\sQ$ and their higher direct images,
so the quantity $\sP$, and it determines the fixed part $F$ and its
intersection numbers. This proves the theorem.

\subsection{The main proof}\label{s!mpf}

The Main Theorem~\ref{th!main} follows formally from the above arguments
together with Proposition~\ref{p!binom}. The plan of the proof: for an
orbifold point $Q$, the local analytic type of $X,\Oh_X(mD)$ is periodic
in $m$, so also the fractional contributions $c_Q(mD)$ of
Theorem~\ref{th!RRD}. The argument of \cite[Theorem 8.5]{YPG} calculates
them by equivariant RR on a global quotient orbi\-fold as Dedekind sums.
Section~\ref{s!DIM} tells us how to replace the Dedekind sums by ice
cream functions, that are integral and Gorenstein symmetric of the given
canonical weight. After subtracting these off, we obtain an integral
valued Hilbert polynomial for $m\gg0$ that is Gorenstein symmetric, to
which Proposition~\ref{p!binom} applies.

\medskip
\paragraph{\bf Step 1}

The local contributions of Theorem~\ref{th!RRD} making up the difference
$\chi(mD)-\RR(mD)$ were calculated in \cite[Theorem 8.5]{YPG} for an
isolated orbifold point of type $\frac1r(a_1,\dots,a_n)$.

\begin{theorem}
\label{thmdis}
Let $X$ be a projective $n$-fold with a basket of isolated cyclic
orbifold points $\sB=\{Q=\frac1r(a_1,\dots,a_n)\}$, and $D$ a
$\Q$-Cartier Weil divisor. Then for $m\in\Z$,
\begin{equation}
\label{eq!cQ}
\chi(X,\Oh_X(mD))=\RR(mD)+\sum_{Q\in\sB}c_Q(mD)),
\end{equation}
where
\begin{equation}
c_Q(mD) = (\si_{r-m}-\si_0)(\tfrac1r(a_1,\dots,a_n)).
\end{equation}
\end{theorem}
Recall the main idea of the proof: by Theorem~\ref{th!RRD}, the
contributions depend only on the analytic type of $(X,mD)$. Thus we can
reduce to the case of a global quotient $X=M/\bmu_r$ having all fixed
points of the same type type $\frac1r(a_1,\dots,a_n)$. The result then
follows by equivariant RR (that is, the Lefschtz fixed point theorem).
\medskip

\paragraph{\bf Step 2}
Ignoring for the moment finitely many initial terms, as traditional in
treating Hilbert polynomials, we replace the genuine Hilbert series
$P_{X,D}(t)=\sum_{m\ge0} h^0(X,mD)t^m$ by the series
$P^\chi_{X,D}(t)=\sum_{m\ge0} \chi(X,mD)t^m$. Since in \eqref{eq!cQ}
$\RR(mD)$ is a polynomial of degree $n$ and the $c_Q(mD)$ are periodic,
summing them gives a term of the form $A(t)/(1-t)^{n+1}$ with
$A(t)\in\Q[t]$ plus periodic terms of the form $B(t)/(1-t^r)$ for each
orbifold point.

Now Corollary~\ref{cor!per} says that the $m$th term in $P\orb(Q,k_X)$
matches the periodic correction $c_Q(mD)$, so that subtracting off our
ice cream functions $P\orb$ reduces us to a formal power series
\begin{equation}
\label{eq!subtract}
P_I^0(t) = P^\chi_{X,D}(t)-\sum_{Q\in\sB}
P\orb\bigl(\tfrac1r(a_1,\dots,a_n),k_X\bigr)
\end{equation}
where $(1-t)^{n+1}P_I^0(t)$ is a polynomial. It follows as usual that
the coefficient of $t^m$ in $P_I^0(t)$ is a polynomial $H(m)$ of degree
$n$ for $m\gg0$, a modified Hilbert polynomial.

\medskip

\paragraph{\bf Step 3}
Now $H(x)$ satisfies the assumptions of Proposition~\ref{p!binom}.
Indeed, it is integer valued because $\chi(\Oh_X(mD))$ and the
coefficients of the power series $P\orb$ are all integers
by~\eqref{PorbInteg!eq}. Moreover, $H(k-x)=(-1)^nH(x)$ because
$\chi(\Oh_X((k-m)D))=(-1)^n\chi(\Oh_X(mD))$ by Serre duality, and we
know by Exercise~\ref{exc!ser} that $\si_{k-m}=(-1)^n \si_m$.

\medskip

\paragraph{\bf Step 4} We define the initial part $P_I$ in terms of the
modified Hilbert polynomial:
\begin{equation}
\label{eq!PI}
P_I(t) = \sum_{m\ge0} H(m)t^m.
\end{equation}
By construction, the two formal power series $P_{X,D}(t)$ and
$P_I(t)+\sum_Q P\orb(t)$ coincide except for an initial segment (since
the first $\rd{\frac c2}$ coefficients of $P\orb(t)$ are zero). This
proves Addendum~\ref{ad}.

\medskip

\paragraph{\bf Step 5} By Appendix~\ref{appendix!Hilb}, $P_I(t)$ has
denominator $(1-t)^{n+1}$ and numerator a palindromic polynomial of
degree $n+k_X+1$, and is therefore determined by its first $\rd{\frac
c2}$ coefficients. Finally, if $R(X,K_X)$ is Gorenstein then these
coefficients are equal to the first $\rd{\frac c2}$ values of
$h^0(X,mD)$. This completes the proof.

\subsection{K3 surfaces and Fano 3-folds}
\label{K3}

Theorem~\ref{th!main} simplifies known results on K3s and Fano 3-folds
(see Alt{\i}nok, Brown and Reid \cite{Si}). Let $(S,D)$ be a polarized
K3 surface with a basket of orbifold points $\sB=\{\frac1r(a,r-a)\}$.
\cite{YPG}, Appendix to Section~8, gives
\begin{equation}
\si_i=\frac{r^2-1}{12r}-\frac{\overline{bi}(r-\overline{bi})}{2r},
\end{equation}
where $ab=1\mod r$ and $\overline{\phantom{\Si}}$ denotes the
smallest nonnegative residue $\mod r$. By Theorem~\ref{th!eurosim},
\[
\InvMod\Bigl((1-t^a)(1-t^{r-a}),\frac{1-t^r}{1-t}\Bigr)\equiv
-\frac1{2r}\sum_{i=1}^{r-1}\overline{bi}(r-\overline{bi})t^i
\]
and
\[
\InvMod\Bigl(\frac{(1-t^a)(1-t^{r-a})}{(1-t)^2},\frac{1-t^r}{1-t}\Bigr)
\equiv
-\frac{(1-t)^2}{2r}\sum_{i=1}^{r-1}\overline{bi}(r-\overline{bi})t^i.
\]

Applying RR for surfaces \cite[Theorem 4.6]{Si} gives the Hilbert series
\begin{equation} \label{eq!K3RR}
P_S(t)=\frac{1+t}{1-t}+\frac{t+t^2}{(1-t)^3}\cdot\frac{D^2}2
-\sum_{\sB} \frac1{1-t^r}\sum_{i=1}^{r-1} \frac{\overline{bi}(r-\overline{bi})}{2r}t^i.
\end{equation}

We can parse $P_S(t)$ into the ice cream functions of
Theorem~\ref{th!main} as follows. Comparing the coefficients of $t$ in
\eqref{eq!K3RR} yields
\begin{equation}
D^2=2g-2+\sum_{\sB}\frac{b(r-b)}r,
\end{equation}
where the genus $g$ is defined by $P_1=h^0(S,\Oh_S(D))=g+1$. Then
$P_S(t)=P_I+\sum_{\sB}P\orb$, where
\begin{equation}
P_I=\frac{1+(g-2)t+(g-2)t^2+t^3}{(1-t)^3}
 = \frac{1+t}{1-t}+(g-1)\frac{t+t^2}{(1-t)^3},
\end{equation}
and one checks as above that
\begin{align}
P\orb &= \frac{\InvMod\Bigl(\frac{(1-t^a)(1-t^{r-a})}{(1-t)^2},
\frac{1-t^r}{1-t},2\Bigr)}
{(1-t)^2(1-t^r)} \\
&= \frac{t+t^2}{(1-t)^3}\cdot \frac{b(r-b)}{2r}
-\frac1{1-t^r}\sum_{i=1}^{r-1}
\frac{\overline{bi}(r-\overline{bi})}{2r}t^i. \label{eq!bibar}
\end{align}
Indeed, the coindex is $c=3$ and the numerator of \eqref{eq!bibar} is
supported in $[2,\dots,r]$.

\begin{cor} Let $V$ be a $\Q$-Fano $3$-fold with basket
$\sB=\left\{\frac1r(1,a,r-a)\right\}$ of terminal quotient
singularities. The Hilbert series of its anticanonical ring is
$P_V(t)=P_I+\sum_{\sB}P\orb$, with
\begin{equation}
P_I=\frac{1+(g-2)t+(g-2)t^2+t^3}{(1-t)^4}
\end{equation}
where $h^0(-K_X)=g+2$ and $-K^3=2g-2+\sum \frac{b(r-b)}r$, and
\begin{equation}
P\orb =\frac{\InvMod\Bigl(\frac{(1-t)(1-t^a)(1-t^{r-a})}{(1-t)^3},
\frac{1-t^r}{1-t},2\Bigr)}
{(1-t)^3(1-t^r)}.
\end{equation}
\end{cor}

\paragraph{\bf Proof} By \cite[Theorem 4.6]{Si} the Hilbert series of
$(V,-K_V)$ equals
\begin{equation}
P_V(t)=\frac{1+t}{(1-t)^2}-\frac{t+t^2}{(1-t)^4}\cdot\frac{K_V^3}2
-\sum_{\sB} \frac1{(1-t)(1-t^r)}\sum_{i=1}^{r-1}
\frac{\overline{bi}(r-\overline{bi})}{2r}t^i.
\end{equation}
The coefficient of $t$ gives the stated value of $-K_V^3$. Clearly Fano
\hbox{$3$-folds} and K3 surfaces have coindex 3 and the same InverseMod
poly\-nomials
$\InvMod\Bigl(\frac{(1-t)(1-t^b)(1-t^{r-b})}{(1-t)^3},
\frac{1-t^r}{1-t}\Bigr)
\equiv\InvMod\Bigl(\frac{(1-t^b)(1-t^{r-b})}{(1-t)^2},
\frac{1-t^r}{1-t}\Bigr)$.
\qed

\begin{exc}\rm Consider the general weighted projective hypersurfaces
\begin{itemize}
\item $S_5\subset\PP(1,1,1,2)$ with an orbifold point of type
$\frac12(1,1)$ at $Q=(0,0,0,1)$;
\item $S_7\subset\PP(1,1,2,3)$ with basket $\{\frac12(1,1),\frac13(1,2)\}$;
\item $S_{11}\subset\PP(1,2,3,5)$ with basket
$\{\frac12(1,1),\frac13(1,2),\frac15(2,3)\}$.
\end{itemize}
All three are K3 surfaces and have $k_{S_i}=0$ and $c=3$. Their Hilbert
series parsed as $P_{S_i}(t)=P_I+\sum_{\sB_i}P\orb$ are as follows
\begin{align*}
P_{S_5}(t)& = \tfrac{1-t^5}{(1-t)^3(1-t^2)}
= \tfrac{1+t^3}{(1-t)^3}+\tfrac{t^2}{(1-t)^2(1-t^2)}, \\
P_{S_7}(t)& = \tfrac{1-t^7}{(1-t)^2(1-t^2)(1-t^3)} \\
& = \tfrac{1-t-t^2+t^3}{(1-t)^3}+\tfrac{t^2}{(1-t)^2(1-t^2)}+\tfrac{t^2+t^3}{(1-t)^2(1-t^3)}, \\
P_{S_{11}}(t)& =\tfrac{1-t^{11}}{(1-t)(1-t^2)(1-t^3)(1-t^5)} = \tfrac{1-2t-2t^2+t^3}{(1-t)^3} \\
& \qquad +\tfrac{t^2}{(1-t)^2(1-t^2)}+\tfrac{t^2+t^3}{(1-t)^2(1-t^3)}+
\tfrac{2t^2+t^3+t^4+2t^5}{(1-t)^2(1-t^5)}.
\end{align*}
\end{exc}

\section{Towards the nonisolated case} \label{s!hd}

This final section speculates on the shape of the Hilbert series of
quasismooth orbifolds with higher dimensional orbifold loci, and
discusses some partial results. We now have abundant experience of
working with these, and definitive results in special cases, such as the
Calabi--Yau 3-fold orbi\-folds of Buckley's thesis \cite{BuSz}. Our
conjectures in the case of 1-dimensional orbifold locus are fairly
specific, and in principle within reach of our methods, although we do
not yet venture formal proofs.

It is traditional to discuss Hilbert functions in terms of Zariski's
notion of {\em cyclic polynomial}, an integral function $H(n)$
represented for $n\gg0$ by $r$ polynomials $f_0,\dots,f_{r-1}$ according
to $n$ modulo~$r$. In the isolated orbifold case discussed so far, the
$f_i$ differ only in their constant terms, so that $H(n)$ behaves
periodically with period $r$, giving the Hilbert series $P(t)=\sum
H(n)t^n$ a simple pole at $\bmu_r$. In the more general case, the $f_i$
differ by terms that grow, giving $P(t)$ higher order poles. The order
of poles corresponds to one plus the dimension of the strata: for if $X$
has a $\frac1s$ orbifold stratum of dimension $\nu$, its graded ring
$R(X,D)$ has at least $\nu+1$ generators $x_i$ of degree $b_i$ divisible
by $s$, and then $P(t)$ usually has a pole of order $\nu+1$ at the
primitive $s$th roots of~1.

{\em Dissident points} are nonisolated orbifold points $Q$ where the
inertia group jumps. Experiments and the results of \cite{BuSz} and
\cite{zh} suggest that for these, the fractional strictly periodic
contribution $\frac\De{1-t^r}$ given by generalized Dedekind sums
(Definition~\ref{df!Dkd}) can be replaced with an integral term, at the
expense of modifying the contributions corresponding to adjacent strata.

\subsection{Calabi--Yau $3$-folds}
The following results illustrate these points. A quasismooth Calabi--Yau
$3$-fold $(X,D)$ has orbi\-fold locus consisting of
\begin{enumerate}
\renewcommand{\labelenumi}{(\alph{enumi})}
\item curves $C$ of generic transverse type $\frac1s(a,s-a)$;
\item points $Q$ of type $\frac1r(a_1,a_2,a_3)$ with
$a_1+a_2+a_3\equiv0\mod r$. If $\hcf(r,a_i)=s_i>1$ then $Q$ is a
dissident point on a $\frac1{s_i}$ curve $C$.
\end{enumerate}

To handle the periodicity of the Hilbert series of $X$, we have the
choice between two alternative strategies. First, a partial fraction
decomposition over $\Q$ with parts having small denominators, directly
linked to the strata of the orbifold locus.

\begin{theorem}\label{th!CY1}
Let $(X,D)$ be as above. Then its Hilbert series is
\begin{equation}
P_{X,D}(t) = \I_{X,D}+\sum_Q \II_Q+\sum_C(\III_C+\IV_C).
\end{equation}
Here the parts are
\begin{equation}
\I_{X,D}=1+\frac t{(1-t)^2}\,\frac{Dc_2}{12}
+\frac{(t+4t^2+t^3)}{(1-t)^4}\,\frac{D^3}6
\end{equation}
with $Dc_2$ and $D^3$ as in RR. (They are, however, rational numbers;
the same applies to the degree $DC$ below.)
\begin{equation}
\II_Q=\frac{\sum(\si_i-\si_0)t^i}{1-t^r}
 =\frac{\De\bigl(\frac1r(a_1,a_2,a_3)\bigr)}{1-t^r},
\end{equation}
with $\si_i=\si_i\bigl(\frac1r(a_1,a_2,a_3)\bigr)$ the Dedekind sum for $Q$
(we note that $\si_0=0$ and $\si_i=-\si_{r-i}$ by Exercise~\ref{exc!ser}).
\begin{equation}
\III_C = \left(\frac{st^s\De}{(1-t^s)^2}
+\frac{t\De'}{1-t^s} - \frac{\si_0 t}{(1-t)^2}\right)DC,
\end{equation}
where $\si_0$ and $\De$ are now the Dedekind sums for the transverse
section $\frac1s(a,s-a)$ of\/ $\Ga$, $\De'=d\De/dt$, and
$\si_0\bigl(\frac1s(a,s-a)\bigr)=\frac{s^2-1}{12s}$.
\begin{equation}
\IV_C= \frac{N_C}{72s\tau_C}\times\frac{B}{1-t^s},
\end{equation}
where the second term has numerator
\begin{multline}
B=t\cdot\InvMod\Bigl(t(1-t^a)^2(1-t^{s-a}),\frac{1-t^s}{1-t}\Bigr) \\
-t\cdot\InvMod\Bigl(t(1-t^a)(1-t^{s-a})^2,\frac{1-t^s}{1-t} \Bigr) \in\Q[t]
\end{multline}
with support $[t,\dots,t^{s-1}]$, determined by
\begin{equation}
(1-t^a)^2(1-t^{s-a})^2B \equiv t^a - t^{s-a} \mod \frac{1-t^s}{1-t},
\end{equation}
and the coefficient $\frac{N_C}{12s^2\tau_C}\in\Q$ depends on the
topology of a tubular neighbor\-hood of $C$ in $X$ (as described in
\cite{BuSz}, Theorem~2.1). The dissident points give rise to the extra
denominator $\tau_C$ and, in spirit, $N_C\in\Z$ is the difference of
degrees of the isotypical components of the normal bundle to $C$;
interchanging $a\bij s-a$ sends $N_C\mapsto -N_C$.
\end{theorem}

This result follows on replacing the individual terms $P_m$ in the
formulas of \cite{BuSz}, Theorem~2.1 by their corresponding Hilbert
series $\sum P_mt^m$ in closed form (with further calculations that are
somewhat involved); it also follows from the results of \cite{zh},
Section~5.1. We omit the details, since our main aim is to motivate the
philosophy of higher-dimensional ice cream, and the detailed statements
are not really the issue.

The second partial fraction decomposition has bigger denominators: each
part has denominator $\prod(1-t^a)$ a product of $n+1=4$ factors. The
parts are one further step removed from the topological characters of
$(X,D)$ appearing in RR. The advantage, however, is that each part is
integral and Gorenstein symmetric of the same degree $k_X=0$.

\begin{theorem}\label{th!CY2}

Let $X,D$ be as above. Then its Hilbert series is
\begin{equation}
P_X(t)=P_I+\sum_QP\orb(Q,k_X) +\sum_C A_C +\sum_CB_C,
\end{equation}
where, as in Theorem~\ref{th!main}, the initial part $P_I$ depends on
the first $\rd{\frac c2}=2$ plurigenera and $P\orb(Q,k_X)$ are ice cream
functions at the point strata;
\begin{equation}
A_C = \frac{P\orb(\tfrac1s(a,s-a),s)}{1-t^s}\, \de C
\end{equation}
with $\de C$ an integer corresponding to degree of $C$ appropriately
modified by the dissident points (see \ref{ss!1dim}), and
\begin{equation}
B_C = \frac{\Num B_C}{(1-t)^3(1-t^s)},
\end{equation}
with numerator an integral palindromic polynomial of symmetric degree
$s+3$ and short support $[3,\dots,s]$.
\end{theorem}

The quantities $\de C$ and $\Num B_C$ are deduced by recombining the
result of Theorem~\ref{th!CY1}, and can be calculated in any particular
example without difficulty, but the formulas are cumbersome to state.
The part $B_C$ depends on the degree $DC$ and of the isotypical
components of its normal bundles, as modified by the dissident points
(see \ref{ss!1dim}). It introduces a periodicity mod~$s$, whereas $A_C$
grows linearly in $s$. We return to the general significance of these
two expressions and the relation between them in \ref{ss!1dim}.

\begin{exa} \rm
The Calabi--Yau 3-fold
\begin{equation}
X_{40}\subset\PP(2,5,8,10,15)_{\Span{x,y,z,t,u}}.
\end{equation}
has the $\frac12(1,1)$ curve
$\Ga_2=X_{40}\cap\PP^2(2,8,10)_{\Span{x,z,t}}$ of degree $\frac12$, and
the $\frac15(2,3)$ curve $C_5=X_{40}\cap\PP^2(5,10,15)_{\Span{y,t,u}}$
of degree $\frac4{15}$ passing through the $\frac1{15}(2,5,8)$ dissident
point $P_u$. One computes the two alternative expressions for its
Hilbert series:
\begin{equation}
\begin{aligned}
P(t) &=\frac{1-t^{40}}{\prod_{a\in\{2,5,8,10,15\}}(1-t^a)} \\
&=\I+\II+\III_{\Ga_2}+\IV_{\Ga_2}+\III_{C_5}+\IV_{C_5} \\
&=P_I+P\orb\bigl(\tfrac1{15}(2,5,8),0\bigr)
 + A_{\Ga_2}+B_{\Ga_2}+A_{C_5}+A_{C_5}.
\end{aligned}
\end{equation}
Here the parts of the first expression are
\begin{equation}
\begin{aligned}
\I &= 1 + \frac{113}{240}\,\frac{t}{(1-t)^2}+ \frac{1}{1800}\,\frac{t+4t^2+t^3}{(1-t)^4}, \quad
\II = \frac{\De(\frac1{15}(2,5,8))}{1-t^{15}}, \\
\III_2 &= -\frac18\, \frac{t+t^3}{(1-t^2)^2}, \quad \IV_2 = 0, \\
\III_5 &= \frac4{15}\,\frac{1 + t^2 + t^3}{(1-t^5)^2}, \quad
\IV_5 = \frac4{25}\,\frac{t^4-\frac83t^3+t^2-t-5}{1-t^5},
\end{aligned}
\end{equation}
where $\De\bigl(\frac1{15}(2,5,8)\bigr)$ is the Dedekind sum polynomial
of Theorem~\ref{th!eurosim}:
\begin{align*}
\De &= t(1-t^5)\InvMod\bigl(t(1-t^5)\cdot(1-t^2)(1-t^5)(1-t^8),1+t^5+t^{10}\bigr) \\
 &=\tfrac19(t + 2t^2 + t^4 - t^5 - 2t^7 + 2t^8 + t^{10} - t^{11} - 2t^{13} - t^{14}).
\end{align*}

The parts of the second expression are
\begin{align*}
P_I &= 1 + \frac{t^2}{(1-t)^4} = \frac{1-4t+7t^2-4t^3+t^4}{(1-t)^4}, \\
P\orb\bigl(\tfrac1{15}(2,5,8),0\bigr) &=
\frac{t^8 - t^9 + t^{10} - t^{11} + t^{12} - t^{13} + t^{14}}{(1-t)^2(1-t^5)(1-t^{15})}, \\
A_{\Ga_2} &=\frac{P\orb(\tfrac12(1,1),2)}{1-t^2}, \quad B_{\Ga_2}=0, \\
A_{C_5} &= \frac{P\orb(\tfrac15(2,3),5)}{1-t^5}, \quad
B_{C_5} = \frac{-3t^3 + 2t^4 - 3t^5}{(1-t)^3(1-t^5)}.
\end{align*}

As a little exercise, we propose the analogous calculations for the
Calabi--Yau 3-fold $X_{80}\subset\PP(3,4,15,20,38)$ (or other cases from
the vast lists of Kreuzer and Skarke).
\end{exa}

\subsection{A general conjecture}
Let $P(t)=\frac {H(t)}{\prod_{i=1}^N (1-t^{b_i})}$ be a rational function
with integral numerator $H(t)\in\Z[t]$ satisfying Gorenstein symmetry
\eqref{eq!Gsym}. For example, $P$ might be the Hilbert series of a Gorenstein
graded ring $R$ of dimension $n+1$ and canonical weight $k_R$ (more generally, a
finite Gorenstein graded module $M$ over a polynomial ring with
$n+1$-dimensional support).

\begin{conj} \label{c!gen}
Under the above assumptions, $P(t)$ has a unique partial fraction
decomposition of the form
\begin{equation}
\label{eq!pars}
P(t) = \sum_A \frac {N_A}{\prod_{a\in A} (1-t^a)}\,.
\end{equation}
The sum runs over sequences $A=\{a_1,\dots, a_{n+1}\}$ consisting of a
main period $r=a_{n+1}$ and some divisors $a_i\mid r$ (some or all of
the $a_i$ may be $1$ or $r$); each $a_i$ divides one of the original
$b_j$, so that a priori only finitely many $A$ occur. The numerator
$N_A$ of each part is an integral polynomial that is symmetric of degree
$k_A=k+\sum_{a\in A} a$, so that the part as a whole has the same
Gorenstein symmetry; moreover $N_A(t)$ is ``of shortest support'', a
minimal residue modulo
\begin{equation}
F_A=\frac{1-t^r}{\hcf\Bigl(1-t^r,\prod_{a\in A, a<r}(1-t^a)\Bigr)}
\end{equation}
(as in \eqref{eq!AF}) supported in an interval of length $<\deg F_A$
centred at $k_A/2$.

To be on the safe side, we could restrict to quasismooth orbifolds.
\end{conj}

If all the $b_i=1$, there is only one part, and the result follows from
Proposition~\ref{p!binom}. We expect the proof to be formal. The idea is
to take account of the poles of $P(t)$ at roots of unity in terms of its
principal parts. The $A$ part should deal with the highest order
principal part of $P(t)$ at primitive $r$th roots of unity, while
possibly modifying the principal parts of higher order poles at
nonprimitive $r$th roots. The parts document the periodicity of an
integral cyclic polynomial, so should have coefficients in $\Z$.

\begin{remark} \rm
Conjecture~\ref{c!gen} delinks Hilbert series and the topological terms
in RR. The conventional narrative is that RR expresses the coherent
cohomology invariants of a projective variety in terms of topological
data. However, in dimension $\ge4$ one does not necessarily aspire to a
detailed understanding of all the Chern numbers in the RR formula. For
example, no-one claims intimate acquaintance with $c_1^2c_2$, $c_1c_3$
and $c_2^2$ in the Todd genus
\begin{equation}
\Td_4=\tfrac1{720}(-c_1^4+4c_1^2c_2+c_1c_3+3c_2^2-c_4).
\end{equation}
In fact, in work with 3-folds, we commonly treat the quantity
$\frac{Dc_2}{12}$ as a basic invariant, deducing its numerical value
from the plurigenera, rather than the other way around. For a canonical
4-fold (say), the plurigenera $P_1,P_2,P_3$ are just integers, and in
our treatment, the initial part
\begin{equation}
P_I=\frac{1+a_1t+a_2t^2+a_3t^3+a_2t^4+a_1t^5+t^6}{(1-t)^5}
\end{equation}
with $a_1=P_1-5$, $a_2=P_2-5P_1+10$, $a_3=P_3-5P_2+10P_1-10$ holds
comparatively few terrors for us, and is arguably a better starting
point than the Chern numbers; for example, they are integers with no
implicit congruences of the type $12 \mid c_1^2+c_2$.

In the same way, even without tying the orbifold parts $P_A$ very
closely to topological invariants (such as the degree of curve orbifold
strata and the isotypical components of their normal bundle), we have
formulas that depend in principle only on a small basket of integers.

\end{remark}

\subsection{The case of curve orbifold locus} \label{ss!1dim}

Let $(X,D)$ be a quasismooth projectively Gorenstein $n$-fold orbifold
of dimension $n\ge2$ with orbifold locus of dimension $\le1$. As before,
its orbifold strata are
\begin{enumerate}
\renewcommand{\labelenumi}{(\alph{enumi})}
\item curves $\Ga$ of transverse type $\frac1{s}(a_1,\dots,a_{n-1})$;

\item points $Q$ of type $\frac1r(a_1,\dots,a_n)$.
\end{enumerate}

Write $s_i=\hcf(a_i,r)$. Dissident points are characterized by having
some $s_i$ a nontrivial factor of $r$, with $1<s_i<r$. The $x_i$-axis is
then pointwise fixed by $\bmu_{s_i}$, so its image is contained in a
$\frac1{s_i}$ orbifold stratum of $X$. Our assumption on the dimension
of the orbifold locus implies that the $s_i$ are pairwise coprime, and
$Q$ is in the closure of orbifold curve strata $\Ga_i$ of transverse
type $\frac1{s_i}(a_1,\dots,\widehat{a_i},\dots,a_n)$.

We summarize the logic of Conjecture~\ref{c!gen} in this case. We treat
the $\frac1s$ orbifold curves (a) by adding contributions of the form
\begin{equation}
\label{eq!cGa}
c_\Ga(t) = \frac{\Num_{D_C}}{(1-t^s)^2(1-t)^{n-1}}
+\frac{\Num_{N_C}}{(1-t^s)(1-t)^n}
\end{equation}
where the numerators are integral, Gorenstein symmetric of the
appropriate degree, and with short support. We expect to see the
$(1-t^s)^2$ in the denominator for the reason outlined at the start of
Chapter~\ref{s!hd}. In the numerators, $D_C$ and $N_C$ refer to
quantities involving the degree of $C$, respectively of the isotypical
components of its normal bundle. Multiplying \eqref{eq!cGa} by
$1-t^{ms}$, corresponding to taking a transverse section by a general
hypersurface in $|msD|$ for some $m$, leaves $\Num_{D_C}$ distinguished
as the numerator of an isolated orbifold point (times $m\times\deg C$).
The $N_C$ term is destroyed by taking a hyperplane section, and cannot
be recovered after so doing.

We deal with points (b) by putting in ice cream of the form
\begin{equation}
\label{eq!genIce}
P\orb(Q,k_X) =
\frac{\InvMod(A,F,\ga)} {\prod_{a\in[s_1,\dots,s_n,r]}(1-t^a)},
\quad\hbox{with}\quad A=\prod\frac{1-t^{a_i}}{1-t^{s_i}},
\end{equation}
with $F$ as in \eqref{eq!AF}, and $\ga$ chosen to make the numerator
Gorenstein symmetric of degree $k_X+r+\sum s_i$.

The contribution $P\orb(Q,k_X)$ is well defined, integral and Gorenstein
symmetric of degree $k$ (see Proposition~\ref{p!np1}). It has the right
periodicity modulo $r$ by an argument similar to
Corollary~\ref{cor!per}.  The curious point, however, is that when
$s_i>1$, it usually contains contributions with denominator
$(1-t^{s_i})^2(1-t)^{n-1}$ and $(1-t^{s_i})(1-t)^n$ that might at first
sight appear to be native to the curves $C_i$ of type
$\frac1{s_i}(a_1,\dots,\widehat{a_i},\dots,a_n)$ through $P$.

A dissident point $Q$ of type $\frac1r$ on an orbifold $\frac1s$ curve
$\Ga$ commonly forces the degree of $\Ga$ and of the isotypical
components of its normal bundle to become fractional with denominator
$r$, thus adding fractional terms into \eqref{eq!cGa}. Attributing
fractional terms with denominator $(1-t^{s_i})^2(1-t)^{n-1}$ and
$(1-t^{s_i})(1-t)^n$ to the dissident point is the same idea as adding a
global fractional term with denominator $(1-t)^{n+1}$ into the local
contribution from an isolated orbifold point, as discussed in
Caution~\ref{ca!gl}.

The meaning of Proposition~\ref{p!np1} is that $P\orb(Q,k_X)$ can be
viewed as obtained from the Dedekind sum term $\frac{\De}{1-t^r}$ by
multiplying top and bottom of the fraction by $\prod(1-t^{s_i})$, then
folding the numerator back into the required interval. Since the
denominator of $P\orb(Q,k_X)$ is $\prod_{a\in[s_1,\dots,s_n,r]}(1-t^a)$,
the difference
\begin{equation}
\frac{\De}{1-t^r} - P\orb(Q,k_X)
\end{equation}
between the Dedekind polynomial and the ice cream function has
$\prod(1-t^{s_i})^2$ in its denominator. This difference has a unique
partial fraction decomposition defined over $\Q$ with parts having
respective denominators
\begin{equation}
(1-t)^{n-2}(1-t^{s_i})^2, \quad (1-t)^{n-1}(1-t^{s_i})
\quad\hbox{and}\quad (1-t)^n,
\end{equation}
and with numerators of short support. Our assertion is that if we parse
the Hilbert series allocating these local parts to the adjacent curves
and to the initial part, we achieve a decomposition with every part
integral and Gorenstein symmetric. Thus using ice cream $P\orb(Q,k_X)$
in place of the more obvious $r$ periodicity contribution
$\frac{\De}{1-t^r}$ effectively cuts the dissident points out of the
curve strata $C$, modifying their degrees and those of the isotypical
components of their normal bundle to be integral.

We believe that for 1-dimensional orbifold locus, the formal statement
and proof of Conjecture~\ref{c!gen} should be within reach of the
strategies of Buckley's thesis \cite{BuSz}. Her proof in the Calabi--Yau
case involves two ingredients: taking cyclic covers in the style of
\cite{YPG} introduces the Dedekind sums at the dissident points. She
deals with the 1-dimensional orbifold locus by resolving singularities
by standard toric resolutions, then calculating exceptional divisors in
the style of our Section~\ref{s!RRD}. Localising around the orbifold
strata more generally is a stacky business, and when we are forced to
wear that hair shirt, we can also hope for progress by combining the
stacky methods of \cite{Edi} and \cite{zh}.

\bigskip

\noindent Anita Buckley, \\
Department of Mathematics, University of Ljubljana, \\
Jadranska 19, 1000 Ljubljana, Slovenia

\noindent {\it e-mail}: Anita.Buckley@fmf.uni-lj.si

\bigskip

\noindent Shengtian Zhou, \\
Mathematics Institute, University of Warwick, \\
Coventry CV4 7AL, England

\noindent {\it e-mail}: Shengtian.Zhou@googlemail.com

\bigskip

\noindent Miles Reid, \\
Mathematics Institute, University of Warwick, \\
Coventry CV4 7AL, England

\noindent {\it e-mail}: Miles.Reid@warwick.ac.uk

\end{document}